\documentclass{proc-l}

\newtheorem{theorem}{Theorem}[section]
\newtheorem{lemma}[theorem]{Lemma}

\newtheorem{prop}[theorem]{Proposition}
\newtheorem{cor}[theorem]{Corollary}

\usepackage{amsfonts}
\usepackage{amsmath}
\usepackage{amsthm}
\usepackage{amssymb}

\newcommand {\theoremref}[1] {Theorem \ref{#1}}
\newcommand {\lemmaref}[1] {Lemma \ref{#1}}
\newcommand {\propref}[1] {Proposition \ref{#1}}

\newcommand {\be} {\begin{enumerate}}
\newcommand {\ee} {\end{enumerate}}

\newcommand {\Ind} {\mathrm{Ind}}
\newcommand {\Res} {\mathrm{Res}}
\newcommand {\End} {\mathrm{End}}
\newcommand {\Aut} {\mathrm{Aut}}

\newcommand {\Z} {\mathbb{Z}}

\newcommand {\C} {\mathbb{C}}
\newcommand {\F} {\mathbb{F}}

\numberwithin{equation}{section}

\begin{document}

\title{Groups With A Character of Large Degree}
\author{Noah Snyder}
\address{Department of Mathematics, University of California, Berkeley, CA, 94720}
\email{nsnyder@math.berkeley.edu}
\thanks{This material is based upon work supported under a National Science Foundation Research Fellowship.}

\subjclass[2000]{20C15}


\begin{abstract}
Let $G$ be a finite group of order $n$ and $V$ a simple $\C[G]$-module of dimension $d$.  For some nonnegative number $e$, we have  $n=d(d+e)$.  If $e$ is small, then the character of $V$ has unusually large degree.  We fix $e$ and attempt to classify such groups.  For $e \leq 3$ we give a complete classification.  For any other fixed $e$ we show that there are only finitely many examples.
\end{abstract}
\maketitle

\section{Introduction}

Throughout this paper we consider a finite group $G$ of order $n$ and a simple $\C[G]$-module $V$ with dimension $d$ and character $\chi$. It is well known  that $d^2 \leq n$ and that $d \mid n$.  Thus, there must exist a nonnegative integer $e$ such that $d(d+e)=n$.  Our approach is to fix $e$ and consider all groups and simple modules such that $d(d+e)=n$.  If $e$ is small relative to $d$ then $\chi$ has large degree.

Our main result is:

\begin{theorem} \label{intro}
Let $G$ be a finite group of order $n$ with a simple $\C[G]$-module $V$ of dimension $d$ and $d(d+e)=n$.
\begin{enumerate}
\item If $e=0$, then $G$ is trivial.
\item If $e=1$ then $G$ is a doubly transitive Frobenius group or a cyclic group with two elements.
\item If $e>1$, then $n \leq ((2e)!)^2$.
\end{enumerate}
\end{theorem}

After finishing this manuscript I discovered that for $e=1$ and $e=2$ this result appears in \cite{Berkovich}, however, that argument does not generalize to larger $e$.

The $e=0$ case is trivial.  The first interesting case is $e=1$, which we treat in section $2$.  In \propref{first-result} we prove part (2) of \theoremref{intro}.  Doubly transitive Frobenius groups were classified by Zassenhaus \cite{Zassenhaus}, so \theoremref{intro} gives a satisfactory classification.

In section $3$ we show that unless $d$ and $e$ satisfy a certain divisibility condition there exists a normal subgroup $N \neq \{1\}$ such that $N$ acts trivially on each simple $\C[G]$-module other than $V$.  In section $3$ we consider the former case and prove \theoremref{intro} in this case.

In section $4$ we turn to the case where such a normal subgroup exists. This was studied by Gagola in \cite{Gagola}, independently by Burtsev in \cite{Burtsev}, and later by Kuisch and van der Waall in \cite{KvdW1}, \cite{Kuisch}, and \cite{KvdW2}.  (For a summary of Burtsev's results in English translation see \cite{Belonogov}.)  We adapt and strengthen a few of their results in \theoremref{maintheorem} to give a purely group theoretic necessary and sufficient condition.  No sufficient condition appeared in the earlier papers.

In section $5$ we use the structure of a certain $p$-group to prove \theoremref{intro} for any fixed $e>1$.  First we prove a weaker version of part $(3)$ of \theoremref{intro} $n \leq e^{2e^2}$, and then we strengthen that result to $n \leq ((2e)!)^2$.   Finally, in section $6$, we give an explicit list of all groups of order $n$ with a character of degree $d$ such that $d(d+2)=n$ or $d(d+3)=n$.

I would like to thank Hendrik Lenstra for many long and fruitful conversations without which this paper would not exist and for his patient comments on drafts of this paper.  I would also like to thank Rob van der Waall and Yakov Berkovich for sharing their knowledge of the relevant literature and Jack Schmidt whose GAP calculations first showed that the $e=3$ case was tractable.

\section{A Preliminary Result}

The first case to consider is $e=0$.  Since the sum of the squares of the dimensions of the simple modules equals the order of the group, $d^2 = n$ implies that there is only one simple module.  This occurs exactly when $G$ is trivial.  

The next case is $e=1$.  This is a particularly nice case since a complete classification is possible.  It is also extraordinary since it is the only $e$ for which there are infinitely many groups.  The results of this section can be derived from the more general results later, but it is helpful to see the outline in this simpler situation first.  (See \cite{Berkovich} for another proof which is elegant, but hard to generalize.)

\begin{prop} \label{first-result}
Let $G$ be a finite group of order $n$.  The following two conditions are equivalent: 
\be
\item There exists a simple $\C[G]$-module $V$ of dimension $d$ with $d(d+1)=n$.

\item $G$ is a semidirect product of a normal subgroup $N$ of order $d+1$ and a group $H$ of order $d$ which acts freely and transitively by conjugation on the set $N-\{1\}$.  Also, unless $d=1$, we have that $V \cong \Ind_N^G V'$ where $V'$ is any nontrivial $1$-dimensional $\C[N]$-module.
\ee
\end{prop}
The groups described in $(2)$, except when $d=1$, are called doubly transitive Frobenius groups.  They were classified by Zassenhaus in \cite{Zassenhaus}.

The proof of \propref{first-result} requires two well-known lemmas.  The first follows immediately from the fact that the number of irreducible representations of a group is the same as the number of conjugacy classes.  The second is \cite[7.2.8i, p. 193]{Robinson}.

\begin{lemma} \label{all-conjugate}
A nontrivial normal subgroup  $N$ of a finite group $G$ acts trivially on each of the nontrivial simple $\C[G]$-modules except for one, if and only if the conjugacy classes of $G$ are $\{1\}$, the set $N-\{1\}$, and the preimages of the nontrivial conjugacy classes of $G/N$.
\end{lemma}

\begin{lemma} \label{elementary-abelian}
If $N$ is a normal subgroup of a finite group $G$ such that all nontrivial elements of $N$ are conjugate in $G$, then $N$ is elementary abelian.
\end{lemma}

\begin{proof}[Proof of \propref{first-result}]
First we show that $(2)$ implies $(1)$.  Let $\chi$ be the character of $V = \Ind_N^G V'$.  Using the explicit formula for the induced character we can compute:
$$\begin{array}{c|c}
g &  \chi(g) \\ \hline
1 & d \\
\in N - \{1\} &-1\\
\notin N &0
\end{array}$$

Therefore, we see that $(\chi, \chi) = (d^2 + d)/n = 1.$  Thus, we have that $V$ is a simple $\C[G]$-module.  Notice that $\chi(1)(\chi(1)+1) = d(d+1) = \#G$.  Thus we see that $(2)$ implies $(1)$.

Now we show that $(1)$ implies $(2)$.  Let $\chi$ be the character of $V$.  If $d=1$, then $n=2$ and $G$ is the cyclic group with two elements.  If $d>1$,  note that $d^2 + d^2 > d^2 + d = n$.  Therefore, $G$ has at most one irreducible character of degree $d$.  Thus all the Galois conjugates of $\chi$ are equal to $\chi$.  Hence, $\chi$ is $\Z$-valued.

By character orthogonality, we see that $$\sum_{g \in G} \chi(g)^2 = n (\chi, \chi) = n = d^2 + d$$ and $$\sum_{g \in G} \chi(g) = n (\chi, 1) = 0.$$  Therefore, we have that $\sum_{g \neq 1} \chi(g)^2 = d$ and $\sum_{g \neq 1} \chi(g) = -d$.  Adding yields, $$\sum_{g \neq 1} \chi(g) (\chi(g) +1) =0.$$  Since all $\chi(g) \in \Z$, this implies that $\chi(g)$ is $0$ or $-1$ for any $g \neq 1$.

Let $S = \{g \in G : \chi(g) = -1\}$ and $T = \{g \in G : \chi(g) = 0\}.$  Now $\sum_{g \neq 1} \chi(g) = -d$ implies that $\#S = d$ and $\#T =d^2-1$.  

Let $W$ be the kernel of the canonical projection $\C[G] \twoheadrightarrow \End_{\C} V$.  Let $\psi$ be the character of $W$.  Notice that $\psi = \rho_G - \chi(1)\chi$ where $\rho_G$ is the character of the regular representation.  Thus we can read off:
$$\begin{array}{c|c}
g &  \psi(g) \\ \hline
1 & d \\
\in S & d\\
\in T & 0
\end{array}$$

Those $g$ for which $\psi(g)=\psi(1)$ are precisely those $g$ which act trivially on $W$.  Thus the map $G \rightarrow \Aut(W)$ has kernel $N := S \cup \{1\}$.  Further notice that $\# N = d+1$.

Since every simple $\C[G]$-module other than $V$ is a quotient of $W$, the subgroup $N$ acts trivially on each of the simple $\C[G]$-modules except for one.  Such simple modules shall be the main focus of the sequel.

Since the index of $N$ is relatively prime to its order, the Schur-Zassenhaus theorem  (cf. \cite[p. 246]{Robinson}) implies that $N$ has a complement $H$.  Since $$d(d+1) = \# G = \# N \# H = (d+1) \# H,$$ we conclude $\# H = d$.  Since there are only $d$ nontrivial elements of $N$ and they are all in the same $H$-orbit under conjugation, $H$ acts freely and transitively on $N -\{1\}$.

By \lemmaref{all-conjugate} all the non-identity elements of $N$ are conjugate in $G$.  So, by \lemmaref{elementary-abelian}, we see that $N$ is abelian.  Take $V'$ to be any $1$-dimensional $\C[N]$-module.  By the other direction of the theorem, $\Ind_N^G V'$ is a simple $\C[G]$-module of dimension $d$, and thus the unique one $V$.
\end{proof}

Thus the existence of a group $G$ of order $n$ with a simple $\C[G]$-module $V$ of dimension $d$ implies that $n$ is of the form $p^k(p^k-1)$ and $d = p^k -1$.  Conversely, if $n = p^k(p^k-1)$ there is at least one example: $G = \F_q \rtimes \F_q^*$ and $V = \Ind_{\F_q}^G V'$ for $V'$ any nontrivial $1$-dimensional module over the group ring of the additive group of $\F_q$.

\section{Reduction to the Two Main Cases}

In the proof of  \propref{first-result}  we found a subgroup $N$ which acts trivially on each simple $\C[G]$-module other than $V$.  This will be our approach to the general case as well.

\begin{theorem}\label{small}
Let $G$ be a finite group of order $n$, and $V$ a simple $\C[G]$-module of dimension $d$ with character $\chi$ and $d(d+e) = n$.  Then at least one of the following is true
\be
\item $d = e$,
\item $d+e \mid (2e-1)!/e$, or
\item there exists some normal subgroup $N \neq \{1\}$ which acts trivially on all simple $\C[G]$-modules other than $V$.
\ee
\end{theorem}

\begin{proof}
If $d < e$ then $d+e \mid (2e-1)!/e$ which is $(2)$.  If $d=e$ we're obviously in case $(1)$.  So we assume that $d>e$.  Hence, we see that  $d^2 + d^2 > d(d+e)$.  So there can be only one irreducible representation of degree $d$.  Thus its character is $\Z$-valued.

Let $W$ be the kernel of the canonical projection $\C[G] \twoheadrightarrow \End_{\C} V$.  Let $\psi$ be the character of $W$.  Let $N$ be the kernel of the map $G \rightarrow \Aut(W)$.  Since every simple $\C[G]$-module other than $V$ is a quotient of $W$, we see that every element of $N$ acts trivially on each simple $\C[G]$-module except $V$. So unless $N = 1$ we have shown $(3)$.

So assume that $G \rightarrow \Aut(W)$ is an injection.  Notice that $\psi(1) = ed$.  For $g \neq 1$, we have that $\psi(g) = - d \chi(g)$.  But, $\psi(g)$ is bounded in absolute value by $\psi(1)=ed$.  Hence, we see $- e \leq \chi(g) \leq e$.  For all $g \neq 1$ we have $\psi(g) \neq \psi(1) = de$, hence $\chi(g) \neq -e$.  Futhermore, since the trivial $\C[G]$-module is a submodule of $W$, thus for any $g \neq 1$, we have that $\psi(g) \neq -de$.  Thus, we also see that $\chi(g) \neq -e$.

Since the product of the characters of two representations is the character of the tensor product, we see that virtual characters form a ring.  By the previous paragraph, the virtual character $\prod_{-e < i < e} (\chi - i)$ vanishes on any $g \neq 1$. Thus, its inner product with the trivial character is $$\frac {1} {d(d+e)} \prod_{1-e \leq i \leq e} (d - i).$$  Canceling the $d$ and using that $d-i \equiv -e-i  \pmod{d+e}$, we see that $d+e \mid (2e-1)!/e$.
\end{proof}

We will call the three cases in \theoremref{small} the \emph{d=e case}, the \emph{factorial case}, and the \emph{normal subgroup case} respectively.  \theoremref{intro} is trivial in the $\emph{d=e}$ case, and follows easily in the factorial case.

\begin{prop} \label{calculation}
If $e$, $d$ and $n$ are integers greater than $1$ such that $n = d(d+e)$ and $d=e$ or $d+e \mid (2e-1)!/e$, then
$n \leq ((2e)!)^2 \leq e^{4e^2}$.
\end{prop}
\begin{proof}
We see immediately that $d \leq d+e \leq (2e)!$.  For $e>1$, we see that $(2e)! \leq (2e)^{2e} \leq (e^2)^{e^2}$.  Thus, $n = d(d+e) \leq ((2e)!)^2 \leq e^{4e^2}$.
\end{proof}

For any $e$ there is an example of the factorial case occurring.  Take $G$ a cyclic group of order $e+1$, and $V$ the trivial $\C[G]$-module.  However, I was unable to find any examples suggesting that the factorial bound is at all tight.

\section{The Normal Subgroup Case}

In this section we study the structure of a group $G$ with a nontrivial normal subgroup $N$ and a simple $\C[G]$-module $V$ such that $N$ acts trivially on each simple $\C[G]$-module other than $V$.

This situation is the main study of \cite{Gagola}.  Gagola showed that the character of $V$ vanishes on all but one nontrivial conjugacy class if and only if there exists a subgroup $N$ of $G$ such that $N$ acts trivially on each simple $\C[G]$-module except $V$.

\begin{theorem}\label{maintheorem}
Let $G$ be a group of order $n$ and $V$ a simple $\C[G]$-module of dimension $d$.  Define $e$ such that $n = d(d+e)$.   Assume that there exists a normal subgroup $N \neq \{1\}$ such that $N$ acts trivially on each simple $\C[G]$-module other than $V$.  Let $x$ be a nontrivial element of $N$ and $C$ be the centralizer of $x$ in $G$.  Then there exist a prime number $p$, a positive integer $k$ and a non-negative integer $m$ such that:

\begin{enumerate} 
\item $N$ is elementary abelian of order $p^k$,
\item $C$ has order $p^k e^2$, and $d = e(p^k-1)$, and $n = e^2 p^k (p^k-1)$,
\item $C$ is a Sylow $p$-subgroup of $G$ and $e = p^m$,
\item if $H$ is any group such that $N \subseteq H \subseteq C$ and $\#H > p^{k+m}$, then $N \subseteq [H,H]$.
\end{enumerate}
\end{theorem}

Parts $(1)-(3)$ and a weaker version of $(4)$ are proved in \cite[Corollary 2.3]{Gagola} and \cite[Corollary 1.4]{KvdW1} using Clifford's theorem.  We give a new proof using different techniques which can be adapted to prove the converse.

\begin{proof}
By \lemmaref{all-conjugate} and \lemmaref{elementary-abelian}, we conclude that $N$ is elementary abelian. Thus for some $p$ and $k$, we see $\#N = p^k$.  This is $(1)$.

Since all but one of the simple $\C[G]$-modules is a pull-back of a simple $\C[G/N]$-moduule, we have $$\C[G] \cong \C[G/N] \oplus \End_{\C}V.$$  Therefore, we see that $d(d+e) = [G:N] + d^2$.  It follows that $[G:N] = de$.  We also have that $$p^k = \#N = \frac {\#G} {[G:N]} = \frac{d}{e} + 1.$$  Therefore, we see that $d = e(p^k -1)$.  Hence, we see that $$n = [G:N] p^k = dep^k = e^2p^k(p^k-1).$$  Because the conjugacy class of $x$ can be identified with $G/C$, we see that $$[G:C] = \#N-1 = \frac d e.$$   Finally, this implies that $$[C:N] = \frac{[G:N]}{[G:C]} = e^2.$$  This completes the proof of $(2)$.

The isomorphism $\C[G] \cong \C[G/N] \oplus \End_\C V$ implies that $\rho_G = \rho_{G/N} + d \chi$ (where $\rho$ denotes the regular character).  This equation lets us read off the values of the character of $V$, which we denote by $\chi$:
$$\begin{array}{c|c}
g &  \chi(g) \\ \hline
1 & d \\
\in N - \{1\} &-e\\
\notin N &0
\end{array}$$

Consider any subgroup $H$ of $G$ and $\varphi$ any degree $1$ character of $H$ which is nontrivial when restricted to $H \cap N$.  We have that 
\begin{align}
\nonumber (\chi, \varphi)_H &= \frac {1} {\#H} \left(d+e - e\sum_{g\in H \cap N} \varphi(g) \right) \\
\nonumber &= \frac {1} {\#H} (d+e - e\#(H \cap N) (\varphi, 1)_{H \cap N}) = \frac{d+e}{\#H}
\end{align}
So we see that \begin{equation} \tag{*} \frac {ep^k} {\#H} = (\chi, \varphi)_H \in \Z.  \end{equation}

We will apply $(*)$ to several different subgroups in order to prove parts $(3)$ and $(4)$. Choose $q$ a prime different from $p$.  Let $Q$ be a $q$-Sylow subgroup of $C$.  Let $H_q = Q \times \langle x \rangle$.  Let $\varphi$ be a degree $1$ character of $H$ which is trivial on $Q$ and nontrivial on $\langle x \rangle$.  Then, $(*)$ tells us that $$\frac {e p^k}{\#H_q} \in \Z.$$  But the $q$ part of $\#H_q$ is the $q$ part of $e^2$.  So $\#H_q$ must be relatively prime to $q$.  Since this is true for all $q \neq p$, we see that $C$ is a $p$-group, and hence that $e=p^m$.   This is $(3)$.

Finally we apply $(*)$ to a group $H$ such that $N \subseteq H \subseteq C$ and $\#H > p^{k+m}$ and any degree $1$ character $\varphi$ of $H$.  Since  $$\frac {p^{k+m}}{\#H} \notin \Z,$$ $(*)$ implies that $\varphi$ must be trivial when restricted to $N$.  Thus, $N \subseteq [H, H]$.  This is part $(4)$.
\end{proof}

\begin{cor}\label{main-cor}
Let $e$ be a positive integer which is not a power of a prime.  Let $G$ be a finite group of order $n$ and $V$ a simple $\C[G]$-module of $V$ of dimension $d$.   Suppose that $d(d+e)=n$, then $n \leq ((2e)!)^2 \leq e^{4e^2}$.
\end{cor}
\begin{proof}
By \theoremref{maintheorem}, the group $G$ does not have a normal subgroup $N \neq \{1\}$ which acts trivially on each simple $\C[G]$-module except for one.  Therefore, by \theoremref{small}, we know that $d+e \mid (2e)!$.  Therefore, by \propref{calculation},  we conclude that $n \leq ((2e)!)^2 \leq e^{4e^2}$.
\end{proof}

Conditions $(1)$-$(4)$ of \theoremref{maintheorem} are also sufficient.  (The reader interested only in the proof of \theoremref{intro} is encouraged to skip to section $5$.)

\begin{theorem} \label{converse}
Let $G$ be a finite group of order $n = p^{k+2m}(p^k-1)$ with $m \in \Z_{\geq 0}$ and $k \in \Z_{>0}$ .  Let $N$ be an elementary abelian normal subgroup of $G$ of order $p^k$.  Let $d = (p^k-1)p^m$.  Let $x$ be a nontrivial element of $N$ and let $C$ be the centralizer of $x$ in $G$.  Suppose that $C$ is 
a Sylow $p$-subgroup of $G$ of order $p^{k+2m}$.  Suppose that if $H$ is any group such that $N \subseteq H \subseteq C$ and $\#H > p^{k+m}$, then $N \subseteq [H,H]$.

Then there exists a simple $\C[G]$-module $V$ such that $\C[G] \cong \C[G/N] \oplus \End_\C V$.  Furthermore, $\dim V = d$.
\end{theorem}

The crucial step in the proof of \theoremref{maintheorem} was that for any virtual character $\chi$ of $G$ and for any subgroup $H$ of $G$ and any character $\varphi$ of $H$ we know that $(\chi, \varphi) \in \Z$. Brauer's theorem gives a converse to this fact.

\begin{prop} \label{brauer} (cf. \cite[p. 82]{Serre}.)
Let $\chi$ be a class function on $G$.  $\chi$ is a virtual character if and only if, for every elementary subgroup $H$ and for every $1$-dimensional character $\varphi$ of $H$, we have $(\chi, \varphi)_H \in \Z$.
\end{prop}

In order to prove \theoremref{converse} we will build a class function $\chi$ and use \propref{brauer} to show that $\chi$ is a character of a simple $\C[G]$-module $V$ satisfying the required properties.

\begin{proof}[Proof of \theoremref{converse}]


Let $\chi$ be the class function such that $d \chi$ is the character of the virtual $\C[G]$-module $\C[G]-\C[G/N]$.  We can read off the values of $\chi$:
$$\begin{array}{c|c}
g &  \chi(g) \\ \hline
1 & d \\
\in N - \{1\} &-p^m\\
\notin N &0
\end{array}$$


If $H$ is an elementary subgroup of $G$ and $\varphi$ is a degree $1$ character of $H$, then we define
\begin{align} \nonumber  f(H, \varphi) := (\chi, \varphi)_H &= \frac 1 {\# H} \left(p^m(p^k-1) -p^m \sum_{g \in H \cap N -\{1\}} \varphi(g)\right) \\ \nonumber &= \frac {p^m} {\# H} \left(p^k - \# (H \cap N) (\varphi, 1)_{H \cap N} \right).\end{align}

Since $C$ is the centralizer of $x$ in $G$, the conjugacy class of $x$ has $(G:C) = p^k -1 = \#N - 1$ elements.  So all nontrivial elements of $N$ are conjugate.  Thus any other nontrivial element $y \in N$ has a Sylow $p$-group of $G$ as its centralizer.  Thus no nontrivial element of $N$ commutes with any element of order prime to $p$.  In particular, if $H$ is an elementary group then $H$ is a $p$-group or $H$ intersects $N$ trivially or both.

Suppose that $H$ is a $p$-group.  Since $f(H, \varphi)$ does not change under conjugation, we can assume that $H \subseteq C$.  Let $H' = H \cdot N$.   Since $N \subseteq H' \subseteq C$, we know that one of $\#H' \mid p^{m+k}$ and $N \subseteq [H', H']$ holds.   On the one hand, if $\# H' \mid p^{m+k}$, then we also have that $\#H \mid p^{m+k}$.  Thus we see that
\begin{align*} f(H, \varphi) &= \frac{p^{m+k}}{\#H} - \frac{p^m \#(H\cap N)}{\#H}(\varphi, 1)_{H \cap N} \\ &=  \frac{p^{m+k}}{\#H} - \frac{p^m \#N}{\#H'}(\varphi, 1)_{H \cap N} \\ &=  \frac{p^{m+k}}{\#H} - \frac{p^{m+k}}{\#H'}(\varphi, 1)_{H \cap N} \in \Z. \end{align*}

On the other hand, suppose that $N \subseteq [H', H']$.  Thus, we see that $$H' = H \cdot N \subseteq H \cdot [H', H'] \subseteq H'.$$ Hence we have that $H' = H \cdot [H', H']$.  By the Burnside basis theorem, \emph{cf.} \cite[p. 135]{Robinson}, any set of elements of $H'$ which generates $H'/[H', H']$ also generates all of $H'$.  Therefore, we see that $H = H'$.  Thus, we know that $N \subseteq [H,H] \subseteq H$.  Since $N \subseteq [H, H]$, all degree $1$ characters of $H$ are trivial on $N$.  Thus, for any $p$-group $H$ we see that $$f(H, \varphi) = \frac {p^m} {\# H} \left(p^k - p^k\right) = 0 \in \Z.$$

Now we suppose that $N \cap H = \{1\}$. Let $H_p$ be a Sylow $p$-subgroup of $H$. Let $r$ be the index  $[H:H_p] $. Thus, 
\begin{align*}  
f(H, \varphi) &= \frac {p^m(p^k-1)} {\# H} = \frac {p^m} {\# H_p} \frac{p^k-1} {r}  \\ 
&= f(H_p, \Res \varphi) \frac{p^k-1}{r}.
\end{align*}
Since the part of $n$ relatively prime to $p$ is $(p^k-1)$, we see that $r \mid p^k-1$.  We already saw that $f(H_p, \varphi) \in \Z$.  Thus we know that $f(H, \varphi) \in \Z$ for any elementary subgroup $H$ and any degree $1$ character $\varphi$ of $H$.  

Therefore, by \propref{brauer}, $\chi$ is a virtual character.  Since $(\chi, \chi) = 1$ and $\chi(1) > 0$ we see that $\chi$ is the character of a simple $\C[G]$-module $V$.  Finally, since $\chi(1)\chi$ is the character of the virtual module $\C[G] - \C[G/N]$ we conclude that $\C[G] \cong \C[G/N] \oplus \End_\C V$.
\end{proof}

Here is one family of groups which satisfy the conditions of \theoremref{converse}.  Let $k$ be a finite field over the prime field $\F_p$, let $W$ be a finite dimensional vector space over $k$, and let $W^\vee$ be its dual space over $k$.  Let $$G = k \times W \times W^\vee \times k^\times$$ with the multiplication $$(a, v, f, c) \cdot (b, w, g, d) = (a + cb + f(w), v + cw, f + g, cd).$$  Let $N = k \subseteq G$.  Here $C$ is $k \times W \times W^\vee$.  Subgroups $H$ such that $N \subseteq H \subseteq C$ are of the form $k \times U$ for $U$ an $\F_p$-subspace of $W \times W^\vee$.  We compute the commutator: $$[(a, v, f, 1), (b, w, g, 1)] = (f(w)-g(v) ,0,0,1).$$  Thus, condition $(4)$ reduces to proving that if $\#U > \# W$, then the antisymmetric $\F_p$-billinear map $U \times U \rightarrow k$ given by $[(v,f), (g,w)] = f(w)-g(v)$ is surjective.  It is sufficient to prove that the $\F_p$-billinear form $[(v,f), (g,w)] = \pi(f(w)-g(v))$ on the $\F_p$-vector space $U$ is nonzero for every $\F_p$-linear surjection $\pi: k \rightarrow \F_p$.  Thus, condition $(4)$ follows from the well-known fact that isotropic subspaces of a symplectic space have dimension at most half the dimension of the space.

\section{The Structure of the $p$-Group $C$}

In order to prove \theoremref{intro} we need to get a bound on $n$ in terms of $e$ in the normal subgroup case.  To do this we will use the Schreier subgroup lemma to study the structure of $C$.

\begin{theorem}[Schreier] cf. \cite[6.1.8]{Robinson} \label{Schreier} 
Let $G$ be a group generated by $n$ elements.  Let $H$ be a subgroup of $G$ of finite index $i$.  Then $N$ is generated by $i(n-1)+1$ elements.
\end{theorem}

\begin{theorem} \label{yay!}
Let  $G$ be a group of order $n$ and $V$ a simple $\C[G]$-module of dimension $d$ such that $n = d(d+e)$.  If $e>1$, then we have that $n \leq e^{4e^2}$.
\end{theorem}
\begin{proof}
In view of \theoremref{small} and \propref{calculation} we need only consider the case where $G$ and $V$ satisfy the conditions of \theoremref{maintheorem}. Thus, we see that $e=p^m$, there exists a $p$-group $H$ with order $p^{k+m+1}$, and $H$ has an elementary abelian subgroup $N$ of order $p^k$, and $N \subseteq [H,H]$.  

Since $N$ is contained in the Frattini subgroup of $H$ and $\#H/N = p^{m+1}$, we see, by the Burnside basis theorem, that $H$ is generated by $m+1$ elements.  Thus, by \theoremref{Schreier}, we see that $N$ is generated by $p^{m+1}m+1$ elements.  Since $N$ is elementary abelian of rank $k$ we see that $k \leq p^{m+1}m+1$.  To get the explicit bound we calculate:

\begin{align*} \log_p(n) &= \log_p\left(p^{2m+k}(p^k-1)\right) \leq 2m+2k \\
&\leq 2m+2(p^{m+1}m + 1) = 2(m+1) + 2p^{m+1}m \leq 4p^{2m}m = \log_p(e^{4e^2}) 
\end{align*}


\end{proof}

Using a result of Gagola's we can strengthen the bound in  \theoremref{yay!}. 

\begin{theorem} \label{centralizer}
Let $N$, $C$ and $e$ be as in \theoremref{maintheorem}, and further assume that $e>1$ and is not a power of $2$.  The centralizer of $N$ in $C$ is strictly larger than $N$.
\end{theorem}
\begin{proof}
This is Theorem 6.2 in \cite{Gagola}.  When $p > 2$, Gagola gave an elegant two paragraph proof on the bottom of page $379$ (the second paragraph together with the paragraph beginning ``Suppose then that $G/N$ is nonsolvable.  If $p>2$..." where the assumption that $G/N$ is nonsolvable is not used.) For $p=2$, Gagola gave a much longer proof of Theorem \ref{centralizer} using the classification of finite simple groups.  As the rest of this paper is elementary, we avoid reliance on the $p=2$ case.
 \end{proof}

\begin{theorem} \label{linear-bound}
Let $e$, $k$, and $m$ be as in \theoremref{maintheorem}.  If $e>1$ and is not a power of $2$, then we have that $k < 2m$.
\end{theorem}
\begin{proof}
By \theoremref{centralizer}, there exists an element $y \in C-N$ such that $N$ is a proper subgroup of the centralizer of $y$ in $G$, denoted $Z_y(G)$.  Let $[y]$ be the conjugacy class of $y$ in $G$.  By \lemmaref{all-conjugate}, we see that $[y]$ is the inverse image of the conjugacy class of the image of $y$ in $G/N$.  In particular we see that $\#N \mid \# [y]$.  Therefore, we know that $$p^k = \#N \mid \#[y] =\frac{\#G}{\#Z_y(G)} \mid \frac{\# G/N}{p} = p^{2m-1}(p^k - 1).$$  Thus, we see that $k < 2m$.
\end{proof}

\theoremref{linear-bound} tightens the bound significantly in the normal subgroup case of \theoremref{yay!}, however, this gives only a marginal improvement in general, since the factorial bound will dominate. 

\begin{theorem}
Let $G$ be a finite group of order $n$ with a simple $\C[G]$-module $V$ of dimension $d$ and $d(d+e)=n$.  Suppose that $e>1$ is not a power of $2$.  Then, $n \leq ((2e)!)^2$.
\end{theorem}
\begin{proof}
In the $d=e$ case or the factorial case, $(d+e) \mid (2e)!$, so this follows immediately.  In the normal subgroup case, since $e$ is not a power of $2$, by \theoremref{linear-bound} we see that, for $e \geq 3$,  $$n = p^{k+2m}(p^k-1) < p^{6m} = e^6 \leq ((2e)!)^2.$$
\end{proof}

Notice that we do not know of any examples of the factorial case with $d>e$.  This leaves open the possibility of much stronger results.  In particular it is still conceivable that $2d^2 > n$ implies the normal subgroup case or that the polynomial bound in \theoremref{linear-bound} holds in general.

\section{Classification for $e = 2$ and $e=3$}

The bounds which we have already proved in \theoremref{small} and \theoremref{linear-bound} are enough to completely classify all groups and all simple modules with $e=2$ or $e=3$.  However, for $e=2$, recall that \theoremref{linear-bound} depends on the classification of finite simple groups.  To avoid this dependence we will use a result of Taussky \cite{Taussky}.

\begin{theorem}\label{2-group}  cf. \cite[Chapter III, Satz 11.9(a)]{Huppert}
Suppose that $C$ is a $2$-group and that $[C,C]$ has index $4$ in $C$.  Then $[C,C]$ is cyclic.
\end{theorem}

\begin{theorem} \label{e=2}
Let $G$ be a finite group of order $n$ with a simple $\C[G]$-module $V$ of dimension $d$ such that $d(d+2) = n$.  Then $G$ is a cyclic group of order $3$, or $G$ is a nonabelian group of order $8$. 
\end{theorem}
\begin{proof}
If $d=1$ then $n=3$.  Thus $d=1$ implies that $G$ is the cyclic group of order $3$.  If $d=2$ then $n=8$.   If $G$ were abelian then it would have no irreducible degree $2$ characters.  Therefore, $d=2$ implies that $G$ is a nonabelian group of order $8$.  Therefore, it is enough to show that $d \leq 2$.

The argument splits up according to the three cases of \theoremref{small}.  Case $(1)$ states that $d=2$.  Case $(2)$ states that $d+2 \mid  3$.  Hence $d = 1$.  Case $(3)$ states that $G$ satisfies the assumptions of \theoremref{maintheorem}.   Thus there exists a nontrivial elementary abelian subgroup $N$ of $G$ and a Sylow $2$-subgroup $C$ of $G$, such that $N$ has index $4$ in $C$, and $N \subseteq [C, C]$.  Since $C/N$ has order $4$ it must be abelian.  Thus, we see that $N = [C,C]$.  So, by \theoremref{2-group}, we see that $N$ is cyclic.  Since it is also nontrivial and elementary abelian it must have order $2$.   But $d = 2 (\#N-1)$, so we see that $d$ must be $2$.
\end{proof}

There is a nice consequence of \theoremref{e=2}.  For $d>1$, we have that $d(d+3) > (d+1)^2$.  Thus, if $G$ is a group and $V$ is a simple $\C[G]$-module such that $\dim V = \lfloor \sqrt{\#G} \rfloor$ then $G$ is cyclic of order at most $3$, or $G$ is nonabelian of order $8$, or $G$ is a doubly transitive Frobenius group.

\begin{theorem} \label{e=3}
Let $G$ be a finite group of order $n$ with a simple $\C[G]$-module $V$ of dimension $d$ such that $d(d+3) = n$.  Then $G$ is one of the two abelian group of order $4$, or $G$ is the dihedral group of order $10$, or $G$ is one of two groups of order $54$ satisfying the conditions of \theoremref{maintheorem}
\end{theorem}
\begin{proof}
Again we split into cases based on \theoremref{small}. The key additional ingredient is \cite[Cor. to Prop. 8.1]{Serre} which states that the degree of any irreducible representation divides the index of any abelian normal subgroup. Case $(1)$ states that $d=3$.  Thus $n = 18$.  Clearly the Sylow $3$-subgroup is abelian, normal, and of index $2$.  Thus $3 \mid 2$, which is a contradiction.

Now consider case $(2)$.  If $d=1$ then $n=4$.  Thus $d=1$ implies that $G$ is an abelian group of order $4$.  If $d=2$ then $n=10$.   If $G$ were abelian then it would have no degree $2$ irreducible characters.  Therefore, $d=2$ implies that $G$ is the nonabelian group of order $10$.  If $d>3$ and $d+3 \mid 5!/3$, it is easy to see that $d$ is prime and that the Sylow $d$-subgroup is normal and cyclic.  Thus, $d \mid d+3$ which contradicts the assumption that $d>3$.

Finally consider case $(3)$.  By \theoremref{linear-bound}, we see that $G$ must be a group of order $54$ satisfying the conditions of \theoremref{maintheorem}.  The fact that there are exactly two such groups is left as an exercise to the reader.
\end{proof}

{}

\end{document}